





%

%
 \font\twelvebf=cmbx12
 \font\twelvett=cmtt12
 \font\twelveit=cmti12
 \font\twelvesl=cmsl12
 \font\twelverm=cmr12		\font\ninerm=cmr9
 \font\twelvei=cmmi12		\font\ninei=cmmi9
 \font\twelvesy=cmsy10 at 12pt	\font\ninesy=cmsy9
 \skewchar\twelvei='177		\skewchar\ninei='177
 \skewchar\seveni='177	 	\skewchar\fivei='177
 \skewchar\twelvesy='60		\skewchar\ninesy='60
 \skewchar\sevensy='60		\skewchar\fivesy='60
%
%

%
 \font\fourteenrm=cmr12 scaled 1200
 \font\seventeenrm=cmr12 scaled 1440
 \font\fourteenbf=cmbx12 scaled 1200
 \font\seventeenbf=cmbx12 scaled 1440
%
%

%
%
%
\font\tenmsb=msbm10
\font\twelvemsb=msbm10 scaled 1200
\newfam\msbfam

%
\font\tensc=cmcsc10
\font\twelvesc=cmcsc10 scaled 1200
\newfam\scfam

%
\def\seventeenpt{\def\rm{\fam0\seventeenrm}%
 \textfont\bffam=\seventeenbf	\def\bf{\fam\bffam\seventeenbf}}
\def\fourteenpt{\def\rm{\fam0\fourteenrm}%
 \textfont\bffam=\fourteenbf	\def\bf{\fam\bffam\fourteenbf}}
\def\twelvept{\def\rm{\fam0\twelverm}%
 \textfont0=\twelverm	\scriptfont0=\ninerm	\scriptscriptfont0=\sevenrm
 \textfont1=\twelvei	\scriptfont1=\ninei	\scriptscriptfont1=\seveni
 \textfont2=\twelvesy	\scriptfont2=\ninesy	\scriptscriptfont2=\sevensy
 \textfont3=\tenex	\scriptfont3=\tenex	\scriptscriptfont3=\tenex
 \textfont\itfam=\twelveit	\def\it{\fam\itfam\twelveit}%
 \textfont\slfam=\twelvesl	\def\sl{\fam\slfam\twelvesl}%
 \textfont\ttfam=\twelvett	\def\tt{\fam\ttfam\twelvett}%
 \scriptfont\bffam=\tenbf 	\scriptscriptfont\bffam=\sevenbf
 \textfont\bffam=\twelvebf	\def\bf{\fam\bffam\twelvebf}%
 \textfont\scfam=\twelvesc	\def\sc{\fam\scfam\twelvesc}%
 \textfont\msbfam=\twelvemsb	
 \baselineskip 14pt%
 \abovedisplayskip 7pt plus 3pt minus 1pt%
 \belowdisplayskip 7pt plus 3pt minus 1pt%
 \abovedisplayshortskip 0pt plus 3pt%
 \belowdisplayshortskip 4pt plus 3pt minus 1pt%
 \parskip 3pt plus 1.5pt
 \setbox\strutbox=\hbox{\vrule height 10pt depth 4pt width 0pt}}
\def\tenpt{\def\rm{\fam0\tenrm}%
 \textfont0=\tenrm	\scriptfont0=\sevenrm	\scriptscriptfont0=\fiverm
 \textfont1=\teni	\scriptfont1=\seveni	\scriptscriptfont1=\fivei
 \textfont2=\tensy	\scriptfont2=\sevensy	\scriptscriptfont2=\fivesy
 \textfont3=\tenex	\scriptfont3=\tenex	\scriptscriptfont3=\tenex
 \textfont\itfam=\tenit		\def\it{\fam\itfam\tenit}%
 \textfont\slfam=\tensl		\def\sl{\fam\slfam\tensl}%
 \textfont\ttfam=\tentt		\def\tt{\fam\ttfam\tentt}%
 \scriptfont\bffam=\sevenbf 	\scriptscriptfont\bffam=\fivebf
 \textfont\bffam=\tenbf		\def\bf{\fam\bffam\tenbf}%
 \textfont\scfam=\tensc		\def\sc{\fam\scfam\tensc}%
 \textfont\msbfam=\tenmsb	
 \baselineskip 12pt%
 \abovedisplayskip 6pt plus 3pt minus 1pt%
 \belowdisplayskip 6pt plus 3pt minus 1pt%
 \abovedisplayshortskip 0pt plus 3pt%
 \belowdisplayshortskip 4pt plus 3pt minus 1pt%
 \parskip 2pt plus 1pt
 \setbox\strutbox=\hbox{\vrule height 8.5pt depth 3.5pt width 0pt}}

%
\def\twelvepoint{%
 \def\small{\tenpt\rm}%
 \def\normal{\twelvept\rm}%
 \def\large{\fourteenpt\rm}%
 \def\huge{\seventeenpt\rm}%
 \footline{\hss\twelverm\folio\hss}%
 \normal}
\def\tenpoint{%
 \def\small{\tenpt\rm}%
 \def\normal{\tenpt\rm}%
 \def\large{\twelvept\rm}%
 \def\huge{\fourteenpt\rm}%
 \footline{\hss\tenrm\folio\hss}%
 \normal}

\tenpoint

%

%
\catcode`\@=11
%
%
\def\footnote#1{\edef\@sf{\spacefactor\the\spacefactor}#1\@sf
 \insert\footins\bgroup\small
 \interlinepenalty100	\let\par=\endgraf
 \leftskip=0pt		\rightskip=0pt
 \splittopskip=10pt plus 1pt minus 1pt	\floatingpenalty=20000
 \smallskip\item{#1}\bgroup\strut\aftergroup\@foot\let\next}
%
%
%
%
\def\hexnumber@#1{\ifcase#1 0\or 1\or 2\or 3\or 4\or 5\or 6\or 7\or 8\or
 9\or A\or B\or C\or D\or E\or F\fi}
\edef\msbfam@{\hexnumber@\msbfam}
\def\Bbb#1{\fam\msbfam\relax#1}
%
%
%
\catcode`\@=12

\twelvepoint
\scriptfont\msbfam=\tenmsb

%
\font\twelvebf=cmbx12	
\font\ninebf=cmbx9	\font\sevenbf=cmbx7	\font\fivebf=cmbx5

\font\twelvebfit=cmbxti10 at 12pt	
\newfam\bfitfam

\font\twelvebm=cmmib10 at 12pt		\font\tenbm=cmmib10
\font\ninebm=cmmib9	\font\sevenbm=cmmib7	\font\fivebm=cmmib5

\skewchar\twelvebm='177	\skewchar\tenbm='177
\skewchar\ninebm='177	\skewchar\sevenbm='177	\skewchar\fivebm='177

\def\twelvepointbold{\def\bold{
\textfont0=\twelvebf	\scriptfont0=\ninebf	\scriptscriptfont0=\sevenbf
\textfont1=\twelvebm	\scriptfont1=\ninebm	\scriptscriptfont1=\sevenbm
\textfont\bffam=\twelvebf	\textfont\bfitfam=\twelvebfit
\def\rm{\fam\bffam\twelvebf}%
\def\it{\fam\bfitfam\twelvebfit}%
\rm}}

%
\twelvepointbold

\font\twelvemsa=msam10 scaled 1200
\newfam\msafam
\textfont\msafam=\twelvemsa
\def\blacksquare{\fam\msafam\mathchar"7004}

\newcount\EQNO      \EQNO=0
\newcount\FIGNO     \FIGNO=0
\newcount\REFNO     \REFNO=0
\newcount\SECNO     \SECNO=0
\newcount\SUBSECNO  \SUBSECNO=0
\newcount\FOOTNO    \FOOTNO=0
\newbox\FIGBOX      \setbox\FIGBOX=\vbox{}
\newbox\REFBOX      \setbox\REFBOX=\vbox{}
\newbox\RefBoxOne   \setbox\RefBoxOne=\vbox{}

\expandafter\ifx\csname normal\endcsname\relax\def\normal{\null}\fi

\def\Eqno{\global\advance\EQNO by 1 \eqno(\the\EQNO)%
    \gdef\label##1{\xdef##1{\nobreak(\the\EQNO)}}}
\def\Fig#1{\global\advance\FIGNO by 1 Figure~\the\FIGNO%
    \global\setbox\FIGBOX=\vbox{\unvcopy\FIGBOX
      \narrower\smallskip\item{\bf Figure \the\FIGNO~~}#1}}
\def\Ref#1{\global\advance\REFNO by 1 \nobreak[\the\REFNO]%
    \global\setbox\REFBOX=\vbox{\unvcopy\REFBOX\normal
      \smallskip\item{\the\REFNO .~}#1}%
    \gdef\label##1{\xdef##1{\nobreak[\the\REFNO]}}}
\def\Section#1{\SUBSECNO=0\advance\SECNO by 1
    \bigskip\leftline{\bf \the\SECNO .\ #1}\nobreak}
\def\Subsection#1{\advance\SUBSECNO by 1
    \medskip\leftline{\bf \ifcase\SUBSECNO\or
    a\or b\or c\or d\or e\or f\or g\or h\or i\or j\or k\or l\or m\or n\fi
    )\ #1}\nobreak}
\def\Footnote#1{\global\advance\FOOTNO by 1 
    \footnote{\nobreak$\>\!{}^{\the\FOOTNO}\>\!$}{#1}}
\def\SameFootnote{$\>\!{}^{\the\FOOTNO}\>\!$}

\def\References{\bigskip\centerline{\bf REFERENCES}
                \smallskip\copy\REFBOX}
\def\NewRefPage{\setbox\RefBoxOne=\vbox{\unvcopy\REFBOX}%
		\setbox\REFBOX=\vbox{}%
		\def\References{\bigskip\centerline{\bf REFERENCES}
                		\nobreak\smallskip\nobreak\copy\RefBoxOne
				\vfill\eject
				\smallskip\copy\REFBOX}%
		\def\NewRefPage{}}


\newcount\LEMMA	  \LEMMA=0
\newcount\PROP	  \PROP=0
\newcount\THM     \THM=0
\def\Lemma#1#2{\global\advance\LEMMA by 1\smallskip\goodbreak
    {\narrower\narrower\narrower\item{\bf Lemma~\the\LEMMA${}_#1$:}
    {\it #2}\smallskip}%
    \gdef\label##1{\xdef##1{\nobreak Lemma~\the\LEMMA${}_#1$}}}
\def\Proposition#1{\global\advance\PROP by 1\smallskip\goodbreak
    {\narrower\narrower\narrower\item{\bf Proposition~\the\PROP:}
    {\it #1}\smallskip}}
\def\Theorem#1#2{\global\advance\THM by 1\smallskip\goodbreak
    {\narrower\narrower\narrower\item{\bf Theorem~\the\THM${}_#1$:}
    {\it #2}\smallskip}%
    \gdef\label##1{\xdef##1{\nobreak Theorem~\the\THM${}_#1$}}}
\def\Proof#1{{\narrower\narrower\narrower\item{Proof:}
    {#1}\nobreak\hfill$\blacksquare$\smallskip}}
\def\Heading#1#2{\smallskip\goodbreak
    {\narrower\narrower\narrower\item{\bf #1:}
    {\it #2}\smallskip}}
\def\MultiRef#1{\global\advance\REFNO by 1 \nobreak\the\REFNO%
    \global\setbox\REFBOX=\vbox{\unvcopy\REFBOX\normal
      \smallskip\item{\the\REFNO .~}#1}%
    \gdef\label##1{\xdef##1{\nobreak[\the\REFNO]}}}
\def\Eqalignno{\global\advance\EQNO by 1 &(\the\EQNO)%
    \gdef\label##1{\xdef##1{\nobreak(\the\EQNO)}}}
\def\NoRef#1{\global\advance\REFNO by 1%
    \global\setbox\REFBOX=\vbox{\unvcopy\REFBOX\normal
      \smallskip\item{\the\REFNO .~}#1}%
    \gdef\label##1{\xdef##1{\nobreak[\the\REFNO]}}}

\input epsf
\def\Fig#1#2#3#4#5{\global\advance\FIGNO by 1 Figure~\the\FIGNO#5
    \topinsert
    \centerline{\epsfysize=#4\epsffile[#3]{#2}}
    {\bigskip\hsize=5.5in\hskip\parindent
     \vbox{\small\item{\bf Figure~\the\FIGNO:}{#1}}}
    \bigskip\endinsert}

\def\RR{{\Bbb R}}
\def\CC{{\Bbb C}}
\def\HH{{\Bbb H}}
\def\OO{{\Bbb O}}

\def\nxn{$n \times n$}
\def\bar{\overline}
\def\Tr{{\rm tr\,}}
\def\Re{{\rm Re}}

\def\rA{r[A]}
\def\rAbar{r[\bar{A}]}


\rightline{22 July 1998}
\bigskip

\null\bigskip
\centerline{\large\bf THE OCTONIONIC EIGENVALUE PROBLEM}
\bigskip

\centerline{Tevian Dray}
\centerline{\it Department of Mathematics, Oregon State University,
		Corvallis, OR  97331, USA}
\centerline{\tt tevian{\rm @}math.orst.edu}
\medskip
\centerline{Corinne A. Manogue}
\centerline{\it Department of Physics, Oregon State University,
		Corvallis, OR  97331, USA}
\centerline{\tt corinne{\rm @}physics.orst.edu}

\bigskip\bigskip
\centerline{\bf ABSTRACT}
\midinsert
\narrower\narrower\noindent
We discuss the eigenvalue problem for $2\times2$ and $3\times3$ octonionic
Hermitian matrices.  In both cases, we give the general solution for real
eigenvalues, and we show there are also solutions with non-real eigenvalues.
\endinsert
\bigskip

\Section{INTRODUCTION}

Finding the eigenvalues and eigenvectors of a given matrix is one of the basic
techniques in linear algebra, with countless applications.  The simplest case
is that of (complex) Hermitian matrices, generalizing the familiar case of
(real) symmetric matrices.  This simple case is nevertheless very important,
for instance in quantum mechanics, where the fact that such matrices have real
eigenvalues allows them to represent physically observable quantities.

The eigenvalue problem is usually formulated over a field, typically either
the real numbers $\RR$ or the complex numbers $\CC$.  We consider here the
generalization to the other normed division algebras, namely the quaternions
$\HH$ and the octonions $\OO$.  We find that most of the basic properties are
retained, provided they are reinterpreted to take into account the lack of
commutativity of $\HH$ and $\OO$, and the lack of associativity of $\OO$.

Our main result is the solution of the real eigenvalue problem for $3\times3$
octonionic Hermitian matrices, also known as Jordan matrices.  It is
straightforward to show
\Ref{H. H. Goldstine \& L. P. Horwitz,
{\it On a Hilbert Space with Nonassociative Scalars},
Proc.\ Nat.\ Aca.\ {\bf 48}, 1134 (1962).} \label\Horwitz
that such matrices admit 24 real eigenvalues, corresponding to eigenvectors
which are independent over~$\RR$.  We show that these eigenvalues do not
satisfy the characteristic equation even though the matrix itself does.
Instead they generically come in 6 sets of multiplicity 4 rather than the
expected 3 sets of multiplicity 8
\Footnote{After this work was completed, we discovered the existence of 
earlier work (in Russian) by Ogievetski\u{\i}~%
\Ref{O. V. Ogievetski\u{\i},
{\it A Characteristic Equation for $3\times3$ Matrices over the Octonions},
Uspekhi Mat.\ Nauk {\bf 36}, 197--198 (1981);
reviewed in Mathematical Reviews 83e:15017.} \label\OVO
which also claims 6 real eigenvalues for such matrices.}%
.  We further show how to generalize the notion of orthogonality to the
nonassociative case, recovering the standard decomposition of a Hermitian
matrix in terms of its eigenvalues and eigenvectors.

We begin in Section~2 with a review of the standard eigenvalue problem for
real and complex Hermitian matrices, and then consider the quaternionic
eigenvalue problem in Section~3.  A brief discussion of the properties of
octonions and octonionic matrices appears in Section~4, after which the
octonionic eigenvalue problem is considered for $2\times2$ and $3\times3$
octonionic Hermitian matrices in Section~5 and Section~6, respectively.
Finally, we discuss our results in Section~7.

\Section{THE STANDARD EIGENVALUE PROBLEM}

The eigenvalue problem as usually stated is to find solutions $\lambda,v$ to
the equation
$$A v = \lambda v \Eqno$$\label\Orig
for a given square matrix $A$.  The basic properties of the eigenvalue problem
for \nxn\ complex Hermitian matrices are well-understood.
\Footnote{We could just as well have started with the case of real symmetric
matrices.}

\Lemma\CC{An \nxn\ complex Hermitian matrix $A$ has $n$ real eigenvalues
(counting multiplicity).}

\Proof{We give here only the proof that the eigenvalues are real.  Let $A$,
$v$, $\lambda$ satisfy \Orig, with $A^\dagger=A$.  Then
$$\bar\lambda v^\dagger v
   = (A v)^\dagger v
   = v^\dagger A v
   = \lambda v^\dagger v
  \Eqno$$\label\Real
so that if $v\ne0$ we have $v^\dagger v\ne0$, which forces
$\bar\lambda=\lambda$.}

\Lemma\CC{Eigenvectors of an \nxn\ complex Hermitian matrix $A$ corresponding
to different eigenvalues are orthogonal.}

\Proof{For $m=1,2$, let $v_m$ be an eigenvector of $A=A^\dagger$ with
eigenvalue $\lambda_m$.  By the previous lemma, $\lambda_m\in\RR$.  Then
$$\lambda_1 v_1^\dagger v_2
  = (A v_1)^\dagger v_2
  = v_1^\dagger A v_2
  = \lambda_2 v_1^\dagger v_2
  \Eqno$$\label\OrthoT
Then either $\lambda_1=\lambda_2$ or $v_1^\dagger v_2=0$.}

\Lemma\CC{For any \nxn\ complex Hermitian matrix $A$, there exists an
orthonormal basis of $\CC^n$ consisting of eigenvectors of $A$.}

\Proof{If all eigenvalues have multiplicity one, the result follows from the
previous lemma.  The Gram-Schmidt orthogonalization process can be used on any
eigenspace corresponding to an eigenvalue with multiplicity greater than one.}

\noindent
These lemmas are equivalent to the standard result that a complex Hermitian
matrix can always be diagonalized by a unitary transformation.  It is
important for what follows to realize that the form of the proofs given above
relies on both the commutativity and the associativity of $\CC$.

Combining the above results, it is easy to see that any (complex) Hermitian
matrix $A$ admits a decomposition in terms of an orthonormal basis
of eigenvectors.

\Theorem\CC{Let $A$ be an \nxn\ complex Hermitian matrix.  Then $A$ can be
expanded as
$$A = \sum_{m=1}^n \lambda_m v_m v_m^\dagger \Eqno$$\label\Decomp
where $\{v_m; ~m=1,...,n\}$ is an orthonormal basis of eigenvectors
corresponding to eigenvalues $\lambda_m$.}

\Proof{By the previous lemma, there exists an orthonormal basis $\{v_m\}$ of
eigenvectors.  It then suffices to check that
$$\sum_{m=1}^n \lambda_m v_m v_m^\dagger v_k = \lambda_k v_k \Eqno$$
But this follows by direct computation using orthonormality.}

\noindent
Furthermore, the set of eigenvalues $\{\lambda_m\}$ is unique, and the (unit)
eigenvectors are unique up to unitary transformations in the separate
eigenspaces (which reduce to multiplication by a complex phase for eigenvalues
of multiplicity one).

\Section{THE QUATERNIONIC EIGENVALUE PROBLEM}
\LEMMA=0
\THM=0

The quaternions $\HH$ double the dimension of the complex numbers by adding
two additional square roots of $-1$, usually denoted $j$ and $k$.  The
multiplication table follows from
$$i^2=j^2=k^2=-1 \qquad ij=k=-ji$$
and associativity; note that $\HH$ is not commutative.  Equivalently, $\HH$
can be viewed via the Cayley-Dickson process as the sum of 2 copies of the
complex numbers
$$\HH = \CC + \CC j\Eqno$$
with $k$ being defined by $k=ij$.

The eigenvalue problem \Orig\ for Hermitian matrices $A$ over $\Bbb H$
immediately yields the first unexpected result: The eigenvalues need not be
real.  An example is given by
$$\pmatrix{~~1&i\cr -i&1} \pmatrix{1\cr k}
  = \pmatrix{1-j\cr k-i}
  = (1-j) \pmatrix{1\cr k}
  \Eqno$$\label\Unex
Furthermore, because of the lack of commutativity, multiples of eigenvectors
are not necessarily eigenvectors.
For instance, the vector
$$v_1 = \pmatrix{~\sqrt{2}\cr 1-i} \Eqno$$
is an eigenvector of the matrix
$$A_1=\pmatrix{~~0&1+i\cr 1-i&~~0} \Eqno$$
with eigenvalue $\sqrt{2}$, but $j v_1$ is not an eigenvector of $A_1$.  This
example illustrates an important point: We must distinguish between right and
left multiplication.  Since
$$A (vq) = (Av) q \Eqno$$
by associativity, {\it right\/} multiples of eigenvectors are indeed
eigenvectors.  For example, $v_1 j$ is an eigenvector of the matrix $A_1$
above, with the same eigenvalue ($\sqrt{2}$).

Similarly, we must carefully distinguish between the {\it left\/} eigenvalue
problem \Orig\ and the {\it right\/} eigenvalue problem
$$A v = v\lambda \Eqno$$\label\Master
It turns out that all the {\it right\/} eigenvalues of (quaternionic) Hermitian
matrices are real.

\Lemma\HH{The right eigenvalues of an \nxn\ quaternionic Hermitian matrix
$A=A^\dagger$ are real.}

\Proof{This is just a careful rewrite of \Real.  Explicitly,
$$\bar\lambda (v^\dagger v)
   = (\bar\lambda v^\dagger) v
   = (A v)^\dagger v
   = (v^\dagger A) v
   = v^\dagger (A v)
   = v^\dagger (v \lambda)
   = (v^\dagger v) \lambda
  \Eqno$$\label\HReal
and the result follows since $v^\dagger v\in\RR$.}

\noindent
In the above, $\bar\lambda$ is now the quaternionic conjugate of
$\lambda$ and $\dagger$ denotes (quaternionic) Hermitian conjugation.
Once the eigenvalues have been shown to be real, orthogonality
of eigenvectors with different eigenvalues follows as in the complex case.

\Lemma\HH{Right eigenvectors of an \nxn\ quaternionic Hermitian matrix $A$
corresponding to different eigenvalues are orthogonal.}

\Proof{This is just a rewrite of \OrthoT.  Explicitly,
$$\eqalign{
\lambda_1 (v_1^\dagger v_2)
  &= (\lambda_1 v_1^\dagger) v_2
   = (A v_1)^\dagger v_2
   = (v_1^\dagger A) v_2 \cr
  &= v_1^\dagger (A v_2)
   = v_1^\dagger (v_2 \lambda_2)
   = (v_1^\dagger v_2) \lambda_2 \cr
  }\Eqno$$
so that, since $\lambda_m\in\RR$, either $\lambda_1=\lambda_2$ or
$v_1^\dagger v_2=0$.}

\noindent
Finally, since the (right) eigenvalues are real, (right) multiples of
eigenvectors are still eigenvectors.

Putting it all together, we obtain a decomposition of any (quaternionic)
Hermitian matrix $A$ of the form \Decomp, where the real eigenvalues
$\{\lambda_m\}$ and their eigenspaces are again unique.

\Lemma\HH{For any \nxn\ complex Hermitian matrix $A$, there exists an
orthonormal basis of $\HH^n$ consisting of eigenvectors of $A$.}

\Theorem\HH{Let $A$ be an \nxn\ quaternionic Hermitian matrix.  Then $A$
can be expanded as in {\rm\Decomp}, where $\{v_m; ~m=1,...,n\}$ is an
orthonormal basis of eigenvectors of $A$, with real eigenvalues $\lambda_m$.}

\noindent
The proofs of each of these results is identical to the complex case
previously considered.  For eigenvalues of multiplicity one, the (unit)
eigenvectors are unique up to a (right) quaternionic phase.

The right eigenvalue problem over $\HH$ is therefore just a straightforward
extension of the complex eigenvalue problem
[\MultiRef{H. C. Lee,
{\it Eigenvalues and Canonical Forms of Matrices with Quaternion
Coefficients},
Proc.\ Roy.\ Irish Acad.\ {\bf A52}, 253--260 (1949).}%
,\MultiRef{J. L. Brenner,
{\it Matrices of Quaternions},
Pacific J. Math.\ {\bf 1}, 329--335 (1951).}%
,\MultiRef{P. M. Cohn,
{\bf Skew Field Constructions},
Cambridge University Press, 1977.}\label\Cohn	
].  The left eigenvalue problem turns out to be of considerable interest as
well, and will be considered elsewhere in the context of $2\times2$ octonionic
Hermitian matrices
\Ref{Tevian Dray, Jason Janesky, and Corinne A. Manogue,
{\it Octonionic Hermitian Matrices with Non-Real Eigenvalues}
(in preparation).} \label\NonReal
\Footnote{Another argument that the right eigenvalue problem is the natural
one is based on their use in diagonalizing a matrix.  Cohn \Cohn\
considers quaternionic matrices $A$ which are diagonalizable, in the sense
that there exists an invertible matrix $U$ and a diagonal matrix $D$ such that
$$A = U D U^{-1}$$
But this is equivalent to
$$A U = U D$$
so that the columns of $U$ are eigenvectors of $A$ with right eigenvalues
taken from the elements of $D$.  We will return to this issue below.  (Cohn
uses the term {\it left eigenvalue} to describe the eigenvector problem for
row vectors multiplied on the right by $A$; this is completely different from
our use of the same term.)}%
.%

\Section{OCTONIONS AND OCTONIONIC MATRICES}

\Subsection{The Octonions}

The octonions $\OO$ can be viewed via the Cayley-Dickson process as the direct
sum of two copies of the quaternions
\Ref{Richard D. Schafer,
{\bf An Introduction to Nonassociative Algebras},
Academic Press, New York, 1966 \& Dover, Mineola NY, 1995.}\label\Schafer
$$\OO = \HH + \HH\ell \Eqno$$
where $\ell$ is yet another square root of $-1$.  The octonions are thus
spanned by the identity element $1$ and the 7 imaginary units
$\{i,j,k,\ell,i\ell,j\ell,k\ell\}$.  These units can be grouped into
associative ``triples'' in 7 different ways, each of which generates (the
imaginary part of) a quaternionic subspace.  Any three of these imaginary
units which do not lie in a such a triple anti-associate.  The multiplication
table can be neatly summarized by appropriately labeling the 7-point
projective plane, as shown in
\Fig{The representation of the octonionic multiplication table using the
7-point projective plane.  Each of the 7 oriented lines represents a
quaternionic triple.}
{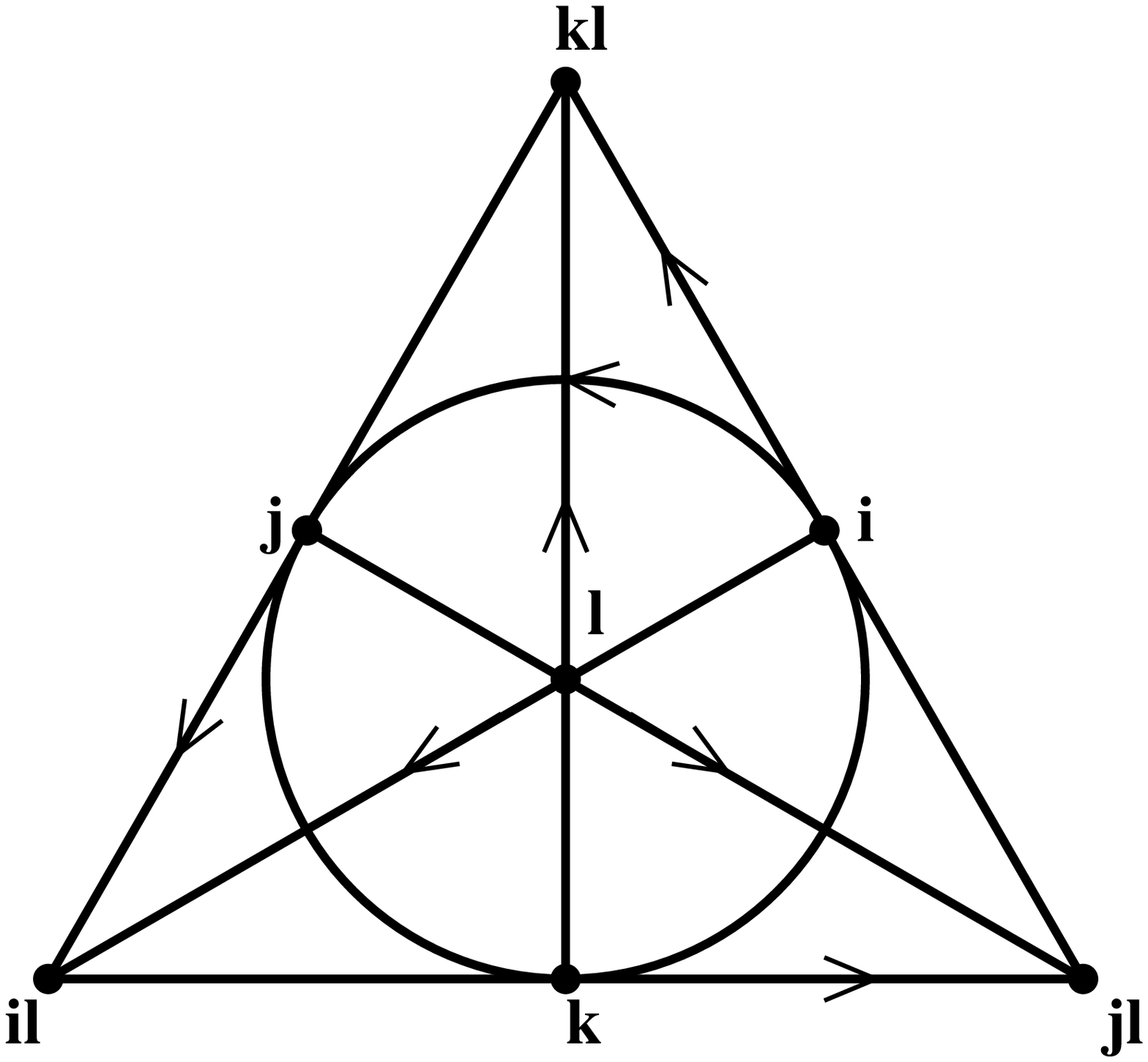}{68 168 543 614}{3in}
{.  For a good introduction to the octonions, including some applications to
physics, see
[\MultiRef{S. Okubo,
{\bf Introduction to Octonion and Other Non-Associative Algebras in Physics},
Cambridge University Press, Cambridge, 1995.}%
,\MultiRef{Feza G\"ursey and Chia-Hsiung Tze,
{\bf On the Role of Division, Jordan and Related Algebras in Particle
Physics},
World Scientific, Singapore, 1996.}]%
.}

When working with small numbers of octonions, it is important to realize that
simplifications take place by virtue of the automorphism group $G_2$ of $\OO$.
For instance, a single octonion may be assumed without loss of generality to
be complex, i.e.\ to lie in the complex subspace of $\OO$ spanned by
$\{1,i\}$.  Similarly, a second octonion can be assumed to lie in the subspace
spanned by $\{1,i,j\}$, and a third can be assumed to lie in the subspace
spanned by $\{1,i,j,k,\ell\}$.  Only when four or more octonions are involved
is it necessary to consider ``generic'' octonions, involving all the basis
directions.

The octonions are not associative.  Nevertheless, since any 2 octonions lie in
a quaternionic subspace, products involving only 2 different octonions (and
their octonionic conjugates) do associate.  For example,
$$p(pq) = p^2 q \Eqno$$
which is a weak form of associativity known as {\it alternativity}.

The squared norm of an octonion $a$ is given by
$$|a|^2 := a\bar{a} \Eqno$$
where $\bar{a}$ denotes the (octonionic) conjugate of $a$.
The commutator of $a$ and $b$ is given by
$$[a,b] := ab - ba \Eqno$$
the associator of $a$, $b$, $c$ is given by
$$[a,b,c] := (ab)c - a(bc) \Eqno$$
and we use $A^\dagger$ to denote the (octonionic) Hermitian conjugate of the
matrix $A$.  Both the commutator and the associator are purely imaginary,
totally antisymmetric, and change sign if any one of their arguments is
replaced by its conjugate.  Another octonionic product with the latter two
properties is given by
\Footnote{$\Phi(a,b,c)$ is in fact the same as the {\it associative 3-form}
[\MultiRef{F. Reese Harvey,
{\bf Spinors and Calibrations}, Academic Press, Boston, 1990.}\label\Reese
,\MultiRef{G. B. Gureirch,
{\bf Foundations of the Theory of Algebraic Invariants},
P. Noordhoff, Groningen, 1964.}]
$$\Phi(a,b,c)
  = \Re (a \times b \times c)
  = {1\over2} \Re \Big( a(\bar{b}c) - c(\bar{b}a) \Big)
  $$
which reduces to the vector triple product when $a$, $b$, $c$ are imaginary
quaternions.
\hfill\break
(Note that $a\times b\times c := {1\over2}\big( a(\bar{b}c)-c(\bar{b}a) \big)$
is the triple cross product, {\it not} the iterated cross product.)}
$$\Phi(a,b,c) = {1\over2} \, \Re( [a,\bar{b}] c) \Eqno$$\label\PhiEq

\Subsection{Octonionic matrices}

The lack of associativity complicates the treatment of matrices with
octonionic entries.  While matrix multiplication can be defined for matrices
of arbitrary size in the usual way, only the $2\times2$ and $3\times3$
octonionic Hermitian matrices form Jordan algebras
\Ref{Nathan Jacobson,
{\bf Structure and Representations of Jordan Algebras},
Amer.\  Math.\  Soc.\ Colloq.\ Publ.\ {\bf 39}, 
American Mathematical Society, Providence, 1968.}\label\Jacobson
.  We therefore limit ourselves to these two cases.

Any complex number $x=x_1 + i x_2$ can be written as the real matrix
$\pmatrix{x_1& -x_2\cr x_2& x_1\cr}$, which allows (complex) matrix
multiplication to be rewritten in terms of real matrices.  This process can be
generalized to the quaternions, but fails for the octonions --- as it must,
since octonionic multiplication is not associative.  But the complex number
$x$ can also be written as the real vector $\pmatrix{x_1\cr x_2}$, leading to
a representation of (complex) multiplication as the product of a (real) matrix
with a (real) vector, and this latter process does generalize to the
octonions.

A similar process can be used to represent matrices over one of these division
algebra as matrices over any smaller such algebra.  Under this transformation
a complex Hermitian matrix is mapped to a symmetric real matrix, and
quaternionic and octonionic Hermitian matrices can be transformed into either
real symmetric matrices or complex Hermitian matrices.  For example, an \nxn\
octonionic Hermitian matrix can be mapped to an $8n\times8n$ symmetric real
matrix (or a $4n\times4n$ complex Hermitian matrix).  It would seem as if we
could therefore reduce the eigenvalue problem to the real case, but this is
misleading for several reasons.

First of all, while the real formalism guarantees the existence of real
eigenvalues, it does {\it not\/} rule out the possibility that eigenvalues
might exist which are not real.  This is because the general eigenvalue
problem over a normed division algebra transforms into a {\it matrix\/}
equation, which only reduces to an ordinary eigenvalue problem for real
eigenvalues.  As we will see, octonionic Hermitian matrices admit (right)
eigenvalues which are not real.

Secondly, in the real formalism it is not very easy to determine the
multiplicity of the real eigenvalues.  One might expect $n$ octonionic
eigenvectors with at most $n$ different real eigenvalues.  The real formalism
{\it does\/} guarantee us $8n$ independent (over $\RR$) eigenvectors with real
eigenvalues, but we are {\it not} guaranteed that the $8n$ eigenvalues occur
with multiplicity $8$ (or a multiple thereof).  In fact, we will see below
that this is not the case.

The final drawback of this approach is that the orthogonality of eigenvectors
with different eigenvalues does {\it not\/} follow from the corresponding
statement on the real, transformed eigenvectors.  This is because the
transformation does not preserve the inner product, which is real in one case
and not in the other.  Nevertheless, the real parts of the inner products do
agree, so that
$$v^\dagger w + w^\dagger v = 0 \Eqno$$\label\OrthoR
As we will see below, there is a stronger orthogonality condition on
eigenvectors with different eigenvalues, which generalizes the usual notion of
orthogonality.

For all of these reasons, we choose to work directly with the octonionic
matrices.

\Section{$\bold 2\times2$ OCTONIONIC HERMITIAN MATRICES}
\LEMMA=0
\THM=0

Are the eigenvalues of octonionic Hermitian matrices real?  Consider first the
special case where $A$ and $v$ (but not necessarily $\lambda$) in \Master\ lie
in a quaternionic subspace.  Then \HReal\ still holds, and alternativity
allows us to shift the parentheses and conclude that $\lambda$ is real.

For instance, since the diagonal entries of a $2\times2$ Hermitian matrix $A$
are real, its components always lie in a {\it complex} subspace of $\OO$, and
therefore $A$ clearly possesses eigenvectors which lie in the same complex
subspace, and which have real eigenvalues.  Multiplication (on the right) of
these eigenvectors by an arbitrary octonion $q$ leads to a new eigenvector,
which lies in the quaternionic subspace spanned by the single octonionic
direction in $A$ and the octonionic multiple $q$.  Furthermore, this new
eigenvector has the same real eigenvalue as the original eigenvector, that is
$$Av = v\lambda \Longrightarrow A (vq) = (Av) q = (v\lambda) q = (vq) \lambda
  \Eqno$$
since $A$ and $v$ are complex and $\lambda$ is real.

In general, however, the key use of associativity in the middle of the
derivation of \HReal\ is not allowed.  We are thus led to suspect that there
exist octonionic Hermitian matrices which admit (right) eigenvalues which are
{\it not\/} real.  This turns out to be correct, as is shown by the following
example:
$$\pmatrix{~~1&i\cr -i&1} \pmatrix{j\cr \ell}
  = \pmatrix{j+i\ell\cr \ell-k}
  = \pmatrix{j\cr \ell} (1-k\ell)
  \Eqno$$\label\SpecEx

Further details for the case of octonionic Hermitian matrices whose
eigenvalues are not real will appear elsewhere \NonReal.

\Subsection{The real eigenvalue problem}

We now turn to the case of {\it real\/} eigenvalues.  The general $2\times2$
octonionic Hermitian matrix can be written
$$A = \pmatrix{p& a\cr \bar{a}& m\cr} \Eqno$$
with $p,m\in\RR$ and $a\in\OO$, and satisfies its characteristic equation
$$A^2 - (\Tr A) \, A + (\det A) \, I = 0 \Eqno$$
where $\Tr A$ denotes the trace of $A$, and where there is no difficulty
defining the determinant of $A$ as usual via
$$\det A = pm - |a|^2 \Eqno$$
If $a=0$ the eigenvalue problem is trivial, so we assume $a\ne0$.  If we set
$$v=\pmatrix{x\cr y\cr} \Eqno$$\label\vform
then \Master\ leads to
$$y = {\bar{a}x \over ~~|a|^2} \, (\lambda-p) \qquad\qquad
  x =       {ay \over ~~|a|^2} \, (\lambda-m)
  \Eqno$$
from which it follows (unless $v=0$) that
$$\det(\lambda I - A) = 
  \lambda^2 - (\Tr A) \, \lambda + \det A = 0 \Eqno$$\label\LambdaI
as usual.  Eigenvectors can thus be written in either of the forms
$$v = \pmatrix{|a|^2\cr \bar{a}(\lambda-p)\cr} \xi \qquad\qquad
  v = \pmatrix{a(\lambda-m)\cr |a|^2\cr} \xi
  \Eqno$$\label\RealForms
where $\xi\in\OO$ is arbitrary and where $\lambda$ is either of the 2
solutions of \LambdaI.
\Footnote{Note that $x\ne0\ne y$ if $a\ne0$, and that in this case the 2
solutions of \LambdaI\ are distinct.}
This shows that {\it all\/} eigenvectors of $2\times2$ Hermitian matrices with
real eigenvalues are obtained from the usual complex eigenvectors by (right)
multiplication by an arbitrary octonion.

\Lemma{{\OO_2}}{Let $A$ be a $2\times2$ complex Hermitian matrix.  Then $w$ is
an octonionic eigenvector of $A$ with real eigenvalue $\lambda$ if and only if
$w=v\xi$, where $\xi\in\OO$ is arbitrary and where $v$ is a complex
eigenvector of $A$ with the same eigenvalue.} \label\ComplexLemma

\noindent
Since a $2\times2$ octonionic Hermitian matrix contains only one independent
octonion, we can assume without loss generality that any such matrix is
complex, and thus apply this lemma to it.

\Subsection{Orthogonality and Decompositions}

As already noted, since $A$ lies in a complex subspace of $\OO$, it admits a
complete set of complex eigenvectors with real eigenvalues, which can be used
to obtain the decomposition~\Decomp.  But what about a decomposition in terms
of the general solution of the real eigenvalue problem?

The general solution is given by \RealForms.  Choosing the first form, we
obtain a complete set of eigenvectors by considering both solutions
$\lambda_\pm$ to \LambdaI, obtaining
$$v_\pm = \pmatrix{|a|^2\cr \bar{a}(\lambda_\pm-p)\cr} \xi_\pm \Eqno$$
Are these eigenvectors orthogonal?  Direct computation yields
$$v_+^\dagger v_- = - |a|^2 a \, \Big[ a,\xi_+,\xi_- \Big] \Eqno$$
\label\NOrtho
where we have used
$$a [a,x,y]
  = -a[\bar{a},y,\bar{x}] \equiv [a,\bar{a}y,\bar{x}] 
  = [\bar{x},a,\bar{a}y]
  = (\bar{x}a) (\bar{a}y) - |a|^2 \bar{x}y
  \Eqno$$
or equivalently
$$a \Big( (\bar{a} \bar{x}) y \Big) = (\bar{x} a) (\bar{a} y) \Eqno$$
for any octonions $a$, $x$, $y$.  Thus, the eigenvectors are not necessarily
orthogonal in the traditional sense except for the quaternionic eigenvalue
problem, when the associator automatically vanishes.

At first sight, this apparent lack of orthogonality seems to rule out a
decomposition of the form \Decomp.  However, due to the lack of associativity,
what is needed for \Decomp\ to hold is an appropriately generalized notion of
orthogonality, namely

\Heading{Definition}
{Let $v$ and $w$ be two octonionic vectors.  We will say that $w$ is
orthogonal to $v$ if
$$(v v^\dagger) \, w = 0 \Eqno$$\label\Ortho
The vectors $\{v,w\}$ are orthonormal if in addition
$v^\dagger v = 1 = w^\dagger w$.}

Direct computation shows that the eigenvectors $v_\pm$ above are indeed
mutually orthogonal in this sense, which provides a computational proof of the
following lemma.

\Lemma{{\OO_2}}{If $v$ and $w$ are eigenvectors of the $2\times2$ octonionic
Hermitian matrix $A$ corresponding to different real eigenvalues, then
$v$ and $w$ are mutually orthogonal in the sense of {\rm\Ortho}.}
\label\TwoOrtho

\Proof{From \ComplexLemma, we can write
$$v=\hat{v}\alpha, \quad w=\hat{v}\beta \qquad
  (\hat{v},\hat{w}\in\CC; \alpha,\beta\in\OO)$$
where $\CC\subset\OO$ is the complex subspace containing the elements of $A$.
But then $vv^\dagger=|\alpha|^2\hat{v}\hat{v}^\dagger$, and
$$(vv^\dagger)w = |\alpha|^2 (\hat{v}\hat{v}^\dagger) (\hat{w} \beta)$$
which associates since only 2 octonionic directions are involved.  But
$\hat{v}^\dagger \hat{w}=0$ by the usual properties of complex eigenvectors.}

In order for a decomposition of the form \Decomp\ to exist, we also need a
vector version of alternativity, which in fact holds for octonionic vectors of
any size:
\Footnote{For 2-component vectors, this proposition is just the 3-$\Psi$'s
rule of supersymmetry theory
\Ref{Corinne A.~Manogue \& Anthony Sudbery,
{\it General Solutions of Covariant Superstring Equations of Motion},
Phys.\ Rev.\ {\bf D40}, 4073-4077 (1989);
\hfill\break
\mathchardef\Psi="7109
Tevian Dray and Corinne A. Manogue,
{\it Associators and the 3-$\Psi$'s Rule},
(in preparation).}%
.}

\Proposition{For any octonionic vector $v\in\OO^n$,
$$( v v^\dagger ) \, v = v \, (v^\dagger v) \Eqno$$
}

\noindent
This proposition shows in particular that any normalized vector $v$ is an
eigenvector of the matrix $vv^\dagger$ with eigenvalue $1$, as required by
\Decomp.  We conclude that the decomposition \Decomp\ holds unchanged for real
eigenvalues.  We thus have:

\Theorem{{\OO_2}}{Let $A$ be a $2\times2$ octonionic Hermitian matrix.  Then
$A$ can be expanded as in {\rm\Decomp}, where $\{v_1, v_2\}$ are orthonormal
(as per~{\rm\Ortho}) eigenvectors of $A$ corresponding to the real
eigenvalues~$\lambda_m$.} \label\Sub

\Proof{Provided the real eigenvalues of $A$ are distinct, \TwoOrtho\
guarantees the existence of orthonormal eigenvectors, which are also
eigenvectors of the decomposition \Decomp\ with the same eigenvalues, and the
result follows.  But if $A$ has a repeated eigenvalue, it must be a multiple
of the identity matrix, for which the result holds.}

Using the same technique as in \TwoOrtho, it is straightforward to show that
$$(vv^\dagger)(vv^\dagger) = (v^\dagger v)(vv^\dagger) \Eqno$$ \label\Idem
for any $v\in\OO^2$, and that
$$(vv^\dagger)(ww^\dagger) = 0 \Eqno$$ \label\IdemOrtho
if $v$ and $w$ are eigenvectors of $A$ with distinct real eigenvalues, since
each term in parentheses lies in $\CC$.  The decomposition in the preceding
theorem is thus in terms of orthonormal idempotents $v_i v_i^\dagger$.  We
also have
$${vv^\dagger\over v^\dagger v}+{ww^\dagger\over w^\dagger w}=I
  \Eqno$$ \label\TwoIdent
which could be proved directly using the fact that the left-hand side has
repeated eigenvalue~1.  Furthermore, since by \TwoOrtho\ $A$ and $v$ contain
only 2 octonionic directions,
$$(Av) v^\dagger = A (vv^\dagger) \Eqno$$ \label\Assoc
which leads to the following one-line alternative derivation of \Sub\
$$
A = A \left( \sum_{i=1}^2 v_i v_i^\dagger \right)
  = \sum_{i=1}^2 A (v_i v_i^\dagger)
  = \sum_{i=1}^2 (Av_i) v_i^\dagger
  = \sum_{i=1}^2 \lambda_i v_i v_i^\dagger
  \Eqno$$ \label\Elegant
As we will see below, however, this argument relies heavily on \TwoOrtho,
which fails in the $3\times3$ case.

We have not yet discussed whether the orthonormal eigenvectors in the
preceding theorem constitute a basis of $\OO^2$.  For any orthonormal vectors,
we have the following lemma.

\Lemma{{\OO_2}}{Let $v,w\in\OO^2$ be orthonormal in the sense of \Ortho, and
let $g$ be any vector in $\OO^2$.  Then
$$g = (vv^\dagger) \, g + (ww^\dagger) \, g$$
}

\Proof{This follows immediately from \TwoIdent.}

\noindent
In the associative case, this lemma shows how to write any vector $g$ in terms
of its components along $v$ and $w$, thus establishing $\{v,w\}$ as a basis.
One could adopt similar language in the nonassociative case, although the
``component'' of $g$ ``along'' $v$ would no longer point in the $v$ direction.
Nevertheless, this terminology is extremely attractive, as it allows the
Gram-Schmidt orthogonalization process to be used to determine the component
of one vector orthogonal to another.

\Proposition{Let $v,w\in\OO^2$.  Then
$$(vv^\dagger) \left( w - {(vv^\dagger)\over v^\dagger v} \, w \right) = 0
  \Eqno$$ \label\GS
}

\Proof{This follows from the alternativity of $2\times2$ octonionic Hermitian
matrices and \Idem.}

\Section{\bold $3\times3$ OCTONIONIC HERMITIAN MATRICES}
\LEMMA=0
\THM=0

We now turn to the $3\times3$ case.  It is not immediately obvious that
$3\times3$ octonionic Hermitian matrices have a well-defined determinant, let
alone a characteristic equation.  We therefore first review some of the 
properties of these matrices before turning to the eigenvalue problem.  As in
the $2\times2$ case, over the octonions there will be solutions of the
eigenvalue problem with eigenvalues which are not real; we consider here only
the real eigenvalue problem.

\Subsection{Jordan matrices}

The $3\times3$ octonionic Hermitian matrices, henceforth referred to as
{\it Jordan matrices}, form the exceptional Jordan algebra (also called the
Albert algebra) under the Jordan product
\Footnote{The $2\times2$ octonionic Hermitian matrices form a special Jordan
algebra since they are alternative \Jacobson.}
$$A \circ B := {1\over2} (AB + BA) \Eqno$$
which is commutative, but not associative.  A special case of this is
$$A^2 \equiv A \circ A \Eqno$$
and we {\it define}
$$A^3 := A^2 \circ A = A \circ A^2 \Eqno$$

Remarkably, with these definitions, Jordan matrices satisfy the usual
characteristic equation \Reese
$$A^3 - (\Tr A) \, A^2 + \sigma(A) \, A - (\det A) \, I = 0 \Eqno$$\label\Char
where $\sigma(A)$ is defined by
$$\sigma(A) := {1\over2} \left( (\Tr A)^2 - \Tr (A^2) \right) \Eqno$$
and where the determinant of $A$ is defined abstractly in terms of the
Freudenthal product.~%
\Footnote{
The Freudenthal product of two Jordan matrices $A$ and $B$ is given by
\Ref{P. Jordan, J. von Neumann, and E. Wigner,
Ann.\ Math.\ {\bf 36}, 29 (1934);
\hfill\break
H. Freudenthal, Adv.\ Math.\ {\bf 1}, 145 (1964).}
$$A*B = A \circ B - {1\over2} \Big(A\,\Tr(B)+B\,\Tr(A)\Big)
		+ {1\over2} \Big(\Tr(A)\,\Tr(B)-\Tr(A\circ B)\Big)$$
The determinant can then be defined as
$$\det(A) = {1\over3} \, \Tr \Big( (A*A) \circ A \Big)$$
}
Concretely, if
$$A = \pmatrix{p& a& \bar{b}\cr \bar{a}& m& c\cr b& \bar{c}& n\cr} \Eqno$$
\label\Three
with $p,m,n\in\RR$ and $a,b,c\in\OO$ then
$$\eqalign{
\Tr A &= p + m + n \cr
\sigma(A) &= pm + pn + mn - |a|^2 - |b|^2 - |c|^2 \cr
\det A &= pmn + b(ac) + \bar{b(ac)} - n|a|^2 - m|b|^2 - p|c|^2 \cr
  }\Eqno$$\label\ThreeEq

\Subsection{The real eigenvalue problem}

As discussed above, \nxn\ Hermitian matrices over any of the normed division
algebras can be rewritten as symmetric $kn\times kn$ real matrices, where $k$
denotes the dimension of the underlying division algebra, it is clear that a
$3\times3$ octonionic Hermitian matrix must have $8\times3=24$ real
eigenvalues \Horwitz.  However, as we now show, instead of having (a maximum
of) 3 distinct real eigenvalues, each with multiplicity 8, we show that there
are (a maximum of) 6 distinct real eigenvalues, each with multiplicity 4.

The reason for this is that, somewhat surprisingly, a (real) eigenvalue
$\lambda$ of a Jordan matrix $A$ does {\it not\/} in general satisfy the
characteristic equation \Char.
\Footnote{Ogievetski\u{\i} \OVO\ constructed a $6$th order polynomial
satisfied by the real eigenvalues, which he called the characteristic
equation.  This polynomial is presumably equivalent to the modified
characteristic equations (for both values of $r$) given below.}
To see this, consider the eigenvalue equation \Orig, with $A$ as in \Three,
$\lambda\in\RR$, and where
$$v=\pmatrix{x\cr y\cr z\cr} \Eqno$$
Explicit computation yields
$$\eqalignno{
  (\lambda-p) x &= ay + \bar{b}z \Eqalignno\cr \label\EI
  (\lambda-m) y &= cz + \bar{a}x \Eqalignno\cr \label\EII
  (\lambda-n) z &= bx + \bar{c}y \Eqalignno\cr \label\EIII
  }$$
so that
$$(\lambda-p) (\lambda-m) y
  = (\lambda-p)(cz + \bar{a}x)
  = (\lambda-p) cz + \bar{a} (ay + \bar{b}z)
  \Eqno$$
which implies
$$\left[ (\lambda-p) (\lambda-m) - |a|^2 \right] y
  = \bar{a} (\bar{b}z) + (\lambda-p) cz
  \Eqno$$ \label\EIV
Assume first that $\lambda\ne p$.  Using \EI\ and \EIV\ in \EIII\ leads to
$$\eqalign{
\left[ (\lambda-p) (\lambda-m) - |a|^2 \right] & (\lambda-p) (\lambda-n) z \cr
&\hskip-5em
  = \left[ (\lambda-p) (\lambda-m) - |a|^2 \right] (\lambda-p) (bx+\bar{c}y) \cr
&\hskip-5em
  = \left[ (\lambda-p) (\lambda-m) - |a|^2 \right] b (ay+\bar{b}z)
    + (\lambda-p)\bar{c}\! \left[ \bar{a}(\bar{b}z) + (\lambda-p) cz \right] \cr
&\hskip-5em
  = b \left[a \left(\bar{a} (\bar{b}z) + (\lambda-p) cz \right) \right]
    + \left[ (\lambda-p) (\lambda-m) - |a|^2 \right] b (\bar{b}z) \cr
   &+ (\lambda-p)\bar{c}\! \left[ \bar{a}(\bar{b}z) + (\lambda-p) cz \right] \cr
&\hskip-5em
  = (\lambda-p) \bigg[ 
	(\lambda-m) |b|^2 z + (\lambda-p) |c|^2 z
	+ b \Big( a (cz) \Big) + \bar{c} \left( \bar{a} (\bar{b}z) \right)
    \bigg] \cr
  }\Eqno$$
Expanding this out and comparing with \ThreeEq\ results finally in
\Footnote{We have recently been able to generalize this to the case where
$\lambda$ is not real \NonReal.}
$$\eqalign{
\Big[ \det(\lambda I - A) \Big] z 
  &\equiv \left[ 
	\lambda^3 - (\Tr A) \, \lambda^2 + \sigma(A) \, \lambda - \det A
     \right] z \cr
  &= b \Big( a (cz) \Big) + \bar{c} \left( \bar{a} (\bar{b}z) \right)
     - \left[ b(ac) + (\bar{c}\,\bar{a})\bar{b} \right] z \cr
  }\Eqno$$\label\CharLam
Now consider the case $\lambda=p$.  We still have \EIV, which here takes the
form
$$-|a|^2 y = \bar{a} (\bar{b}z) \Eqno$$\label\FI
Inserting this into \EII, we can solve for $x$, obtaining
$$-|a|^2 x = a (cz) + (p-m) \bar{b}z \Eqno$$\label\FII
Finally, inserting \FI\ and \FII\ in \EIII\ yields
$$-\left( |a|^2 (p-n) + |b|^2 (p-m) \right) z
  = b \Big( a (cz) \Big) + \bar{c} \left( \bar{a} (\bar{b}z) \right)
  \Eqno$$
Comparing with \ThreeEq\ and using $\lambda=p$, we see that \CharLam\ still
holds, and thus holds in general.

If $a$, $b$, $c$, and $z$ associate, the RHS of \CharLam\ vanishes, and
$\lambda$ does indeed satisfy the characteristic equation \Char; this will
not happen in general.  However, since the LHS of \CharLam\ is a real multiple
of $z$, this must also be true of the RHS, so that
$$b \Big( a (cz) \Big) + \bar{c} \left( \bar{a} (\bar{b}z) \right)
     - \left[ b(ac) + (\bar{c}\,\bar{a})\bar{b} \right] z
  = rz \qquad\qquad r\in\RR
  \Eqno$$\label\RReal
which can be solved to yield a quadratic equation for $r$ as well as
constraints on $z$.

\Lemma{{\OO_3}}{The real eigenvalues of the $3\times3$ octonionic Hermitian
matrix $A$ satisfy the modified characteristic equation
$$\det(\lambda I - A)
  = \lambda^3 - (\Tr A) \, \lambda^2 + \sigma(A) \, \lambda - \det A
  = r
  \Eqno$$\label\Lameq
where $r$ is either of the two roots of
$$r^2 + 4\Phi(a,b,c) \, r - \Big| [a,b,c] \Big|^2 = 0 \Eqno$$\label\Req
with $a,b,c$ as defined by {\rm\Three} and where $\Phi$ was defined in
{\rm\PhiEq}.}

\NoRef{Tevian Dray and Corinne A. Manogue,
{\it Finding Octonionic Eigenvectors Using {\sl Mathematica}},
Comput.\ Phys.\ Comm.\
(invited paper; submitted).}\label\Find

\Proof{These results were obtained using {\sl Mathematica} to solve \RReal\
by brute force for real $r$ and octonionic $z$ given generic octonions
$a$,~$b$,~$c$ \Find.}

\noindent
Furthermore, provided that $[a,b,c]\ne0$, each of $x$, $y$, and $z$ can be
shown to admit an expansion in terms of 4 real parameters.
\Heading{Corollary 1}
{With $A$ and $r$ as above, and assuming $[a,b,c]\ne0$,
$$z = (\alpha a + \beta b + \gamma c + \delta)
      \left( 1 + {[a,b,c] \, r \over \Big| [a,b,c] \Big|^2} \right)
  \Eqno$$\label\Zeq
with $\alpha,\beta,\gamma,\delta\in\RR$.  Similar expansions hold for $x$ and
$y$.}

\noindent
The real paramaters $\alpha,\beta,\gamma,\delta$ may be freely specified for
one (nonzero) component, say $z$; the remaining components $x,y$ have a
similar form which is then fully determined by \EI--\EIII.

\Heading{Corollary 2}
{The real eigenvalues of $\bar{A}$ are the same as those of $A$.}

\Proof{Direct computation (or \ThreeEq) shows that
$$\det \bar{A} = \det A - 4 \Phi(a,b,c) \Eqno$$
But $-4\Phi(a,b,c)$ is precisely the sum of the roots of \Req, and replacing
$A$ by $\bar{A}$ merely flips the sign of $r$, that is $\rAbar=-\rA$.  Thus,
the 2 possible values of $\det A + \rA$ are precisely the same as those for
$\det\bar{A} + \rAbar$.  Since $\Tr{\bar{A}}=\Tr{A}$ and
$\sigma(\bar{A})=\sigma(A)$, \Lameq\ is unchanged.}

The solutions of \Lameq\ are real, since the corresponding $24\times24$ real
symmetric matrix has 24 real eigenvalues.  We will refer to the 3 real
solutions of \Lameq\ corresponding to a single value of $r$ as a {\it family}
of eigenvalues of $A$.  There are thus 2 families of real eigenvalues, each
corresponding to 4 independent (over $\RR$) eigenvectors.

We note several intriguing properties of these results.  If $A$ is in fact
complex, then the only solution of \Req\ is $r=0$, and we recover the usual
characteristic equation with a unique set of 3 (real) eigenvalues.  If $A$ is
quaternionic, then one solution of \Req\ is $r=0$, leading to the standard set
of 3 real eigenvalues and their corresponding quaternionic eigenvectors.
However, unless $a$, $b$, $c$ involve only two independent imaginary
quaternionic directions (in which case $\Phi(a,b,c)=0=[a,b,c]$), there will
also be a nonzero solution for $r$, leading to a second set of 3 real
eigenvalues.  From the preceding corollary, we see that this second set of
eigenvalues consists precisely of the usual ($r=0$) eigenvalues of $\bar{A}$!
Furthermore, since $A(\ell v)=\ell(\bar{A}v)$ for $A,v\in\HH$ and imaginary
$\ell\in\OO$ orthogonal to $\HH$, the eigenvectors of $A$ corresponding to
$r\ne0$ are precisely $\ell$ times the quaternionic ($r=0$) eigenvectors of
$\bar{A}$.  In this sense, the octonionic eigenvalue problem for quaternionic
$A$ is equivalent to the quaternionic eigenvalue problem for both $A$ and
$\bar{A}$ together.  Finally, if $A$ is octonionic (so that in particular
$[a,b,c]\ne0$), then there are two distinct solutions for $r$, and hence two
different sets of real eigenvalues, with corresponding eigenvectors.  Note
that if $\det{A}=0\ne[a,b,c]$ then all of the eigenvalues of $A$ will be
nonzero!

The final suprise lies with the orthogonality condition for eigenvectors $v,w$
corresponding to different eigenvalues.  It is {\it not\/} true (cf.\ \NOrtho)
that $v^\dagger w=0$, although the real part \OrthoR\ of this expression does
vanish.  However, just as in the $2\times2$ case, what is needed to ensure a
decomposition of the form \Decomp\ is \Ortho, and a lengthy, direct
computation verifies that \Ortho\ holds {\it provided that\/} both
eigenvectors correspond to the same value of $r$.

\Lemma{{\OO_3}}{If $v$ and $w$ are eigenvectors of the $3\times3$ octonionic
Hermitian matrix $A$ corresponding to different real eigenvalues in the same
family (same $r$ value), then $v$ and $w$ are mutually orthogonal in the sense
of {\rm\Ortho}.} \label\MainOrtho

\Proof{The modified characteristic equation \Lameq\ can be used to
eliminate cubic and higher powers of $\lambda$ from any expression.
Furthermore, given two distinct eigenvalues $\lambda_1\ne\lambda_2$,
subtracting the two versions of \Lameq\ and factoring the result leads to the
equation
$$(\lambda_1^2+\lambda_1\lambda_2+\lambda_2^2)
  - \Tr{A} (\lambda_1+\lambda_2) + \sigma(A)
  = 0 \Eqno$$
which can be used to eliminate quadratic terms in one of the eigenvalues.
\Footnote{We used {\sl Mathematica} to implement these simplifications in a
brute force verification of \Ortho\ in this context, which ran for 6 hours on
a SUN Sparc20 with 224 Mb of RAM \Find.}
}

\noindent
For Jordan matrices, we thus obtain {\it two\/} decompositions of the form
\Decomp, corresponding to the two sets of real eigenvalues.  For each, the
eigenvectors are fixed up to orthogonal transformations which preserve the
form \Zeq\ of $z$.  

\Theorem{{\OO_3}}{Let $A$ be a $3\times3$ octonionic Hermitian matrix.  Then
$A$ can be expanded as in {\rm\Decomp}, where $\{v_1, v_2, v_3\}$ are
orthonormal (as per {\rm\Ortho}) eigenvectors of $A$ corresponding to the real
eigenvalues $\lambda_m$, which belong to the same family (same $r$
value).}\label\Main

\Proof{Fix a family of real eigenvalues of $A$ by fixing $r$.  If the
eigenvalues are distinct, then the previous theorem guarantees the existence
of orthonormal eigenvectors, which are also eigenvectors of the decomposition
\Decomp\ with the same eigenvalues, and the result follows.
\hfill\break\indent
If the eigenvalues are the same, the family consists of a single real
eigenvalue $\lambda$ with multiplicity 3.  Then $\Tr(A)=3\lambda$ and
$\sigma(A)=3\lambda^2$.  Writing out these two equations in terms of the
components \Three\ of $A$, and inserting the first into the second, results in
a quadratic equation for $\lambda$; the discriminant $D$ of this equation
satisfies $D\le0$.  But $\lambda$ is assumed to be real, which forces $D=0$,
which in turn forces $A$ to be a multiple of the identity matrix, for which
the result holds.
\hfill\break\indent
The remaining case is when one eigenvalue, say $\mu$, has multiplicity 2 and
one has multiplicity 1.  Letting $v$ be a (normalized) eigenvector with
eigenvalue $\mu$, consider the matrix
$$X = A - \alpha \> v v^\dagger \Eqno$$\label\Two
with $\alpha\in\RR$.  For most values of $\alpha$, $X$ will have 3 distinct
real eigenvalues, whose eigenvectors will be orthogonal by the previous
theorem.  But this means that eigenvectors of $X$ are also eigenvectors of
$A$; the required decomposition of $A$ is obtained from that of $X$ simply by
solving \Two\ for $A$.}

\noindent
Note in particular that for some quaternionic matrices
with determinant equal to zero, one and only one of these two decompositions
will contain the eigenvalue zero.

In the $2\times2$ case, \Idem\ tells us that, for normalized $v$, $vv^\dagger$
squares to itself, and hence is idempotent.  As already noted, the
decomposition of \Sub\ is thus an idempotent decomposisition.  But \Idem\
fails in the $3\times3$ case, so that the decomposition in \Main\ is therefore
{\it not} an idempotent decomposition.

It is nevertheless straightforward to show, in analogy with \TwoIdent, that if
$u$, $v$, and $w$ are orthonormal in the sense of \Ortho, then
$$uu^\dagger + vv^\dagger + ww^\dagger = I \Eqno$$ \label\ThreeIdent
since the left-hand side has eigenvalue $1$ with multiplicity $3$.  This
permits us to view $\{u,v,w\}$ as a basis of $\OO^3$ in the following sense

\Lemma{{\OO_3}}{Let $u,v,w\in\OO^3$ be orthonormal in the sense of \Ortho, and
let $g$ be any vector in $\OO^3$.  Then
$$g = (uu^\dagger) \, g + (vv^\dagger)  \, g + (ww^\dagger) \, g \Eqno$$
}

\Proof{This follows immediately from \ThreeIdent.}

\noindent
However, another consequence of the failure of \Idem\ in the $3\times3$ case
is that the Gram-Schmidt orthogonalization procedure \GS\ no longer works.  It
appears to be fortuitous that we are nevertheless able to find orthonormal
eigenvectors in the $3\times3$ case with repeated eigenvalues; we suspect that
this might fail in general, perhaps already in the $4\times4$ case with an
eigenvalue of multiplicity 3.

\Section{DISCUSSION}

Our main result is \Lameq\ together with \Req, which shows how to modify the
characteristic equation for a $3\times3$ octonionic Hermitian matrix in order
to find its real eigenvalues, and which further shows that there are in
general 2 families of solutions of these equations, each consisting of 3
eigenvalues of multiplicity 4.

We have further shown how to use the corresponding eigenvectors to recover the
usual decomposition \Decomp\ of a Hermitian matrix in terms of its
eigenvectors for both $2\times2$ and $3\times3$ octonionic Hermitian matrices.
In the process we were led to introduce an appropriately generalized notion of
orthogonality, namely \Ortho.

We can relate our notion of orthonormality to the usual one by noting that
$n$ vectors in $\OO^n$ which are orthonormal in the sense \Ortho\ satisfy
$$vv^\dagger + ... + ww^\dagger=I \Eqno$$
If we define a matrix $U$ whose columns are just $v,...,w$, then this
statement is equivalent to
$$U U^\dagger = I \Eqno$$
Over the quaternions, left matrix inverses are the same as right matrix
inverses, and we would also have
$$U^\dagger U = I \Eqno$$
or equivalently
$$v^\dagger v = 1 = ... = w^\dagger w; \qquad v^\dagger w = 0 = ... \Eqno$$
which is just the standard notion of orthogonality.  These two notions of
orthogonality fail to be equivalent over the octonions; we have been led to
view the former as more fundamental.

We can now rewrite the eigenvalue equation \Master\ in the form
$$A U = U D \Eqno$$ \label\Matrix
where $D$ is a diagonal matrix whose entries are the real eigenvalues.
Multiplying \Matrix\ on the left by $U^\dagger$ yields
$$U^\dagger (AU) = U^\dagger (UD) = (U^\dagger U) D \Eqno$$
(since $D$ is real), but this does not lead to a diagonalization of $A$ since,
as noted above, $U^\dagger U$ is not in general equal to the identity matrix.
However, \Main\ can be rewritten as
$$A = U D U^\dagger \Eqno$$
so that in this sense $A$ is diagonalizable.  Furthermore, multiplication of
\Matrix\ on the right by $U^\dagger$ shows that
$$(AU) U^\dagger = (UD) U^\dagger = A = A (UU^\dagger) \Eqno$$
and this assertion of associativity can be taken as a restatement of both
\Sub\ and \Main.  In the $2\times2$ case, this associativity holds for a
single eigenvector $v$, which is \Assoc, and which was used in the one-line
proof \Elegant\ of \Sub.  However, \Assoc\ fails in the $3\times3$ case, and
we are unaware of a correspondingly elegant proof of \Main.

Many of our results were obtained using {\sl Mathematica} \Find.  In fact, the
{\it only} proof we currently have of the $3\times3$ orthogonality result,
namely \MainOrtho, uses {\sl Mathematica} to explicitly perform a horrendous,
but exact, algebraic computation.  While one could hope for a more elegant
mathematical proof of this result, the {\sl Mathematica} computation
nevertheless establishes a result which would otherwise remain for the moment
merely a conjecture.  This is a good example of being able to use the computer
to verify one's intuition when it may not be possible to do so otherwise.
This issue is further discussed in \Find.

Finally, it is intriguing that (some) Hermitian octonionic matrices admit
eigenvalues which are not real.  In particular, this means that octonionic
self-adjoint operators do not necessarily have a (purely) real spectrum.  We
plan to report separately on these matters \NonReal.

\vfill\eject
\bigskip\leftline{\bf ACKNOWLEDGMENTS}\nobreak

We are grateful to Al Agnew for unwittingly suggesting the topic, and to Tom
Craven, Burt Fein, and Tony Sudbery for discussions.


\References

\bye